


\magnification\magstep1
\baselineskip14pt 
\vsize23.5truecm 
\hsize15.0truecm
\hoffset1.0truecm 
\overfullrule=0pt

\font\bigbf=cmbx12 
\font\csc=cmcsc10 

\def\hatt{\widehat}
 
\def\sumin{\sum_{i=1}^n}

\def\d{{\rm d}}
\def\E{{\rm E}} 
\def\Var{{\rm Var}} 

\def\dell{\partial} 
 
\def\arr{\rightarrow} 
\def\half{\hbox{$1\over 2$}}
\def\quart{\hbox{$1\over 4$}}
\def\sixth{\hbox{$1\over 6$}}
\def\section{\medskip}
\def\ET{Efron and Tibshirani}
\def\tr{{\rm t}}
\def\mtrix{\pmatrix} 

\def\ref#1{{\noindent\hangafter=1\hangindent=20pt
  #1\smallskip}}          

\def\cstok#1{\leavevmode\thinspace\hbox{\vrule\vtop{\vbox{\hrule\kern1pt
        \hbox{\vphantom{\tt/}\thinspace{\tt#1}\thinspace}}
        \kern1pt\hrule}\vrule}\thinspace} 
 

\centerline{\bigbf Performance of Efron and Tibshirani's}
\centerline{\bigbf semiparametric density estimator} 

\smallskip
\centerline{\bf Nils Lid Hjort, University of Oslo}
\centerline{\bf December 1995} 

\smallskip
{{\smallskip\narrower\baselineskip12pt\noindent
{\csc Abstract.} 
Recently, \ET\ (1996) proposed a semiparametric density estimator,  
which works by multiplying an initial kernel type estimate 
with a parametric exponential type correction factor,
chosen so as to match certain empirical moments. 
While \ET\ investigate and illustrate many aspects of their method,
the basic questions of performance, and comparison with 
other density estimators, were not directly addressed in their article.  
The purpose of this paper is to provide formulae for 
bias and variance and hence mean squared error for the estimator. 
This additional insight into the method makes it easy to compare
its performance with that of other recently proposed 
semiparametric constructions. 
A brief comparison study is carried out here. 
It indicates that the new method, 
used with lower order polynomials in the exponential correction term, 
is often better than the kernel estimator,
in a reasonable neighbourhood around the normal distribution, 
but that its performance as a density estimator 
is more than equalled by other methods. 

\smallskip\noindent
{\csc Key Words:} \sl 
bias reduction, 
comparative study, 
density estimation, 
multiplicative correction,
semiparametric, 
specially designed exponential families
\smallskip}}

\bigskip
{\bf 1. Introduction and summary.}
Let $X_1,\ldots,X_n$ be independently sampled from a smooth 
density $f$. The traditional kernel method estimator is 
$\hatt f_0(x)=n^{-1}\sumin K_h(x_i-x)$, 
where $K(u)$ is a symmetric probability density kernel function, 
with scaling $K_h(u)=h^{-1}K(h^{-1}u)$ afforded via a smoothing 
parameter $h$. Its basic properties include 
$$\E\hatt f_0\doteq f+\half k_2h^2f''
        \quad {\rm and} \quad
  \Var\hatt f_0\doteq (nh)^{-1}R(K)f-n^{-1}f^2, \eqno(1.1)$$
where $k_2=\int u^2K(u)\,\d u$ and $R(K)=\int K(u)^2\,\d u$;
see for example Scott (1992) or Wand and Jones (1995). 

\ET\ (1996) recently introduced a semiparametric density estimator, 
effectively putting such an initial kernel estimator through an
exponential family. Specifically, consider the parametric model 
$$f(x,\beta)=\hatt f_0(x)\,\hatt c(\beta)^{-1}\exp\{\beta^\tr t(x)\}
        =\hatt f_0(x)\,\hatt c(\beta)^{-1}\exp\{\beta_1t_1(x)
                +\cdots+\beta_pt_p(x)\}, \eqno(1.2)$$
with integration constant 
$$\hatt c(\beta)=\int \hatt f_0(x)\exp\{\beta^\tr t(x)\}\,\d x. $$ 
(If $K$ has bounded support then existence of this and 
similar integrals is guaranteed.) 
This is a specially designed exponential family,
defined in terms of canonical variables $t(x)=(t_1(x),\ldots,t_p(x))^\tr$,
and with $\hatt f_0$ as `carrier density'. 
The \ET\ proposal is to calculate the maximum likelihood estimate
$\hatt\beta$ for the model (1.2), 
while in this calculational step 
ignoring the sampling variability present in the carrier. 
The result is the density estimator 
$$\hatt f(x)=f(x,\hatt\beta)
 =\hatt f_0(x)\,\hatt c(\hatt\beta)^{-1}
        \exp\{\hatt\beta^\tr t(x)\}. \eqno(1.3)$$
Natural choices for $t(x)$ would include lower order polynomials, 
say $(x,x^2)$, $(x,x^2,\allowbreak x^3)$, $(x,x^2,x^3,x^4)$. Using (1.3) 
with $t(x)=(x,x^2,x^3,x^4)$, for example, 
sees to it that the density estimator
has exactly the same mean, variance, skewness and kurtosis as 
the data themselves. 

\ET\ explore and illustrate many aspects of this semiparametric method.
They use discretisation and a link to Poisson regression 
to draw on the large toolkit of generalised linear models to
address issues as choice of $t(x)$,
significance of $\beta_j$s, useful diagnostic tools, and so on. 
Viewed as a smoothing method their work leads in particular
to innovative methods for choosing the bandwidth. 
This framework may also lend itself naturally 
to choosing a locally varying bandwidth, 
using lessons from regression literature; 
see in this connection for example Simonoff (1996, Ch.~6). 
\ET\ also discuss conceptual similarities with other  
recent semiparametric proposals, including the semiparametric
method of Hjort and Glad (1995);
see Sections 2 and 3 below for further comparison. 

One aspect is however conspicuously absent in their article, 
namely that of statistical performance, 
and hence comparison with other methods. 
The main purpose of this paper is to provide 
good approximations to the basic bias and variance quantities, 
and hence to mean squared error, in the traditional 
large-sample framework in which $h\arr0$ while $nh\arr\infty$ 
(these conditions will also be seen to be 
minimum requirements for consistency). 
With formulae for bias and variance we may also 
much more easily compare the actual performance
of the new estimator to others on the market. 

In the following result, which is proved in Section 4, 
a central r\^ole is played by the function 
$$g(x)=(\E t''(X))^\tr\Sigma^{-1}(t(x)-\xi), \eqno(1.4)$$
where $\xi$ and $\Sigma$ are mean vector and variance matrix 
for $t(X)$ under the true density~$f$. 

{\smallskip\sl
{\csc Proposition.} 
Let $f$ be a density with two continuous derivatives 
around a point $x$. 
Suppose the variables $t_j(X_i)^2$ and $t_j''(X_i)$ 
have finite means, and that the kernel function $K$ has bounded support. 
Then the semiparametric estimator (1.3) has 
$$\eqalign{
\E\hatt f(x)&=f(x)+\half k_2h^2\{f''(x)-f(x)g(x)\}+o(h^2)+O(h^2n^{-1}) \cr 
{\sl and\ }\Var\hatt f(x)
        &=(nh)^{-1}R(K)f(x)-n^{-1}f(x)^2+O(hn^{-1}) \cr}\eqno(1.5)$$
as $h\arr0$ and $nh\arr\infty$. 
(If $f$ admits four continuous derivatives around $x$,
then the $o(h^2)$ term is in fact $O(h^4)$.) 
\smallskip}

These results are quite similar to (1.1) for the kernel method itself;
the variance is not affected at all, 
to the order of approximation employed, 
while the bias is of the same order, but different. 
There is indeed a surprising variety of recently proposed estimators,
constructed from different perspectives, with the same property:
the variance is approximately unaffected, and the bias is of 
the form $\half k_2h^2b(x)$ plus smaller order terms. 
Many of these are reviewed and compared in Hjort and Jones (1996b);
see also Sections 2 and 3 below. 
Viewed in this light the \ET\ estimator is yet another 
`$h^2$ bias modification kernel method'. 
Among the immediate consequences is that the rate of convergence cannot 
be better than that of the kernel method itself, or than that of 
the many other methods in this large group of methods; 
the optimal rate for the $h$ bandwidth is the traditional $O(n^{-1/5})$,
leading to optimal mean squared error of size $O(n^{-4/5})$.
Good performance of (1.3) can be expected for densities 
for which $f''-fg$ of (1.5) is small in size.
For the low order polynomial choices of $t(x)$ 
we find in Section 3 that this corresponds to being close to the normal.

It is also important to realise that there can be no uniform improvement;
while the (1.3) method sometimes and perhaps often 
does improve on the kernel estimator,  
there are still many instances where it does not. 
Examples are exhibited in Section 3 where 
this modern and sophisticated estimator will perform worse 
than the more than a century old kernel method (see Peirce, 1873). 

We present competing estimators from recent literature in Section 2,
and summarise a brief comparison study in Section 3.
This study suggests in particular that the new (1.3) estimator
is often outperformed by some of its natural competitors.
Then Section 4 gives the proof of the main result (1.4)--(1.5). 
Finally Section 5 offers some generalisations and supplementing remarks. 

\bigskip
{\bf 2. Similarity to other estimators.} 
We refer to Hjort and Jones (1996b) for further comparisons 
with many other $h^2$ bias modification methods. 
Here we merely single out and describe those 
in this bag of variegated inventions 
that appear to have intentions most similar to those of the \ET\ scheme (1.3).
Estimators $\hatt f_2$, $\hatt f_3$, $\hatt f_4$, $\hatt f_5$ 
described below all have variances of the form $R(K)f(x)/(nh)-f(x)^2/n$ 
plus smaller order terms,
and each has bias of the form $\half k_2h^2b(x)$ plus smaller order terms.
For the present purposes it therefore suffices 
to exhibit the $b(x)$ functions.

The \ET\ method is motivated by matching moments; 
if $X^*$ is drawn from $\hatt f$ of (1.3),
then $t_j(X^*)$ has the same mean as the empirical $\bar t_j$ for the data.
Such moment matching can be achieved in other ways as well.
Jones (1991) discusses a couple of simple ways in which the variance of the
kernel estimate can be matched to the empirical variance $s^2$ of the data,
one of which is to calculate the kernel estimate on 
$x_i+s(s^2+h^2)^{-1/2}(x_i-\bar x)$. 
The resulting estimator, say $\hatt f_2$, 
could be said to have the same aim as the \ET\ method 
with $t(x)=(x,x^2)$. Jones showed  
$$b_2(x)=f''(x)+\sigma^{-2}\{f(x)+(x-\mu)f'(x)\}, \eqno(2.1)$$
with $\mu$ and $\sigma$ being the mean and standard deviation 
parameters (assumed finite) for the density. 

While the (1.3) estimator performs a parametric correction 
on a nonparametric start, Hjort and Glad (1995) 
make a nonparametric correction on a parametric start.
They take a parametric $f(x,\hatt\theta)$, 
say the normal or a mixture of two normals, 
and then multiply with a nonparametric correction factor 
$n^{-1}\sumin K_h(x_i-x)/f(x_i,\hatt\theta)$. 
For this estimator, say $\hatt f_3$, 
$$b_3=f_0(f/f_0)''
        =f''-(f')^2/f+f\{(f'/f-f_0'/f_0)^2
        -f_0''/f_0+(f_0'/f_0)^2\}, \eqno(2.2)$$
where $f_0(x)=f(x,\theta_0)$ is the best parametric approximant to $f(x)$,
see Hjort and Glad (1995, Sections 2 and 3) for the details. 
When maximum likelihood estimators are used, for example,
`best approximation' means minimising the Kullback--Leibler distance
from $f(x)$ to $f(x,\theta)$. 
Different vehicle families give different bias functions $f_0(f/f_0)''$. 

Some of the estimators arrived at 
via local likelihood machinery in Hjort and Jones (1996a)
and via related multiplicative correction methods in Hjort (1996) 
also have appearances similar to (1.3). 
Suppose again that $f(x,\hatt\theta)$ is some initial parametric estimator, 
and consider a local parametric approximation of the form 
$f(t)=f(t,\hatt\theta)r(t,\alpha)$ for $t$ around a given $x$. 
The final density estimator at $x$ is then 
$\hatt f(x)=f(x,\hatt\theta)r(x,\hatt\alpha_x)$, 
where $\hatt\alpha_x$ maximises 
the localised kernel smoothed log-likelihood function
$$\sumin K_h(x_i-x)\log\{f(x_i,\hatt\theta)r(x_i,\alpha)\}
        -n\int K_h(t-x)f(t,\hatt\theta)r(t,\alpha)\,\d t, $$
see Hjort and Jones (1996a), and also Loader (1996a), 
for motivation and particularities of such constructions. 
The simplest case is a local level model
$r(t,\alpha)=\alpha$ for the correction function, which is found to give 
$$\hatt f_4(x)=f(x,\hatt\theta)\hatt\alpha_x
        =\hatt f_0(x){f(x,\hatt\theta)
          \over \int K_h(t-x)f(t,\hatt\theta)\,\d t}. $$
Similarly, if the local parametric correction model uses
both a local level and a local slope, say
$r(t,\alpha_1,\alpha_2)=\alpha_1\exp(\alpha_2(t-x))$, 
then the estimator can be written 
$$\hatt f_5(x)=\hatt f_0(x)
        {f(x,\hatt\theta)\over \int K(u)f(x+hu,\hatt\theta)
          \exp(\hatt\alpha_{2,x}hu)\,\d u}, $$
where $\hatt\alpha_{2,x}$ is the solution to a certain equation.
These estimators have integral $1+O(h^2)$ in general, 
but normalisation is possible via numerical integration 
if deemed important. 

Estimators $\hatt f_4$ and $\hatt f_5$ can be shown 
to have bias factor functions 
$$\eqalign{
b_4&=f''-f(f_0''/f_0), \cr
b_5&=f''-(f')^2/f+f\{(f_0'/f_0)^2-f_0''/f_0\}, \cr} \eqno(2.3)$$
where once more $f_0(x)=f(x,\theta_0)$ is the best approximant 
within the parametric family. 
These results follow from methods in Hjort and Jones (1996a, 1996b).
Again different start families give different estimators
with different bias functions. 
See Section 3 for specific examples. 

The semiparametric estimators above essentially involve 
one nonparametric operation in combination with a globally parametric one. 
We note in passing that it is also possible to 
correct the nonparametric kernel estimator with a 
second nonparametric correction factor.
These could for example look like $\hatt f_3$, $\hatt f_4$ or $\hatt f_5$ 
but with a nonparametric start estimator replacing $f(x,\hatt\theta)$, 
see Jones, Linton and Nielsen (1995) and Hjort (1996). 
Such methods manage to have $h^4$ order bias functions,
while keeping the variance at $O((nh)^{-1})$ level, 
and may outperform the $h^2$ order bias methods,
and in particular the (1.3) estimator, 
at least for large sample sizes and for sufficiently smooth densities. 
Another way of achieving such $O(h^4)$ bias,
and thus superior performance for all large $n$, 
is to use $f(x,\hatt\theta)r(x,\hatt\alpha_x)$ 
as for $\hatt f_4$ and $\hatt f_5$, but with three or four 
local parameters in the $r(t,\alpha)$ term;
see Hjort and Jones (1996a). 

\bigskip
{\bf 3. Specific applications, and comparisons.}
In this section we study performance of the \ET\ method 
with the most natural low order polynomial choices for $t(x)$,
and make comparisons to the similarly-inspired estimators 
$\hatt f_2$, $\hatt f_3$, $\hatt f_4$, $\hatt f_5$
described in Section 2 above. 

\section
{\csc 3.1. Low order polynomial versions.} 
With $t(x)=x$ the $g(x)$ function of (1.4) is zero, 
which means that the performance of the (1.3) estimator is 
identical to that of the kernel estimator.
This is not even an approximate statement,
since inspection shows that $\hatt\beta=0$ in this situation. 

Consider then $(x,x^2)$, for which $t''(x)=(0,2)$. 
Hence, with $\mu$ and $\sigma$ denoting 
mean and standard deviation for $f$, 
$g$ of (1.4) takes the form 
$$g_2(x)=(0,2)\,\Sigma^{-1}\mtrix{x-\mu \cr x^2-(\sigma^2+\mu^2)} 
        =2\sigma^{21}(x-\mu)+2\sigma^{22}\{x^2-(\sigma^2+\mu^2)\}, $$
where the $\sigma^{jk}$s are the elements of $\Sigma^{-1}$,
and where $\Sigma$ is the variance matrix of $(X_i,X_i^2)^\tr$. 
One may express $\sigma^{21}$ and $\sigma^{22}$ in terms of 
the skewness and kurtosis of $f$. 
Observe that the bias factor $f''-fg_2$ vanishes under normality. 
Thus particularly good performance can be expected 
in a nonparametric neighbourhood around the normal distribution. 

Next let $t(x)=(x,x^2,x^3)$, with $t''(x)=(0,2,6x)$. Here 
$$\eqalign{
g_3(x)&=(0,2,6\mu)\,\Sigma^{-1}
\mtrix{x-\mu \cr x^2-(\sigma^2+\mu^2) \cr x^3-\E X^3 \cr} \cr
&=2[\sigma^{21}(x-\mu)+\sigma^{22}(x^2-\E X^2)+\sigma^{23}(x^3-\E X^3)] \cr
&\quad
+6\mu[\sigma^{31}(x-\mu)+\sigma^{32}(x^2-\E X^2)
        +\sigma^{33}(x^3-\E X^3)], \cr}$$
with $\sigma^{jk}$s on this occasion denoting elements of 
the inverse of the variance matrix $\Sigma$ for $(X_i,X_i^2,X_i^3)^\tr$,
and therefore involving the first six moments of $f$. 
The case of $(x,x^2,x^3,x^4)$ may be analysed in the same way,
featuring $\E t''(X)=(0,2,6\mu,12(\mu^2+\sigma^2))^\tr$,
and giving an explicit $g_4$ function, 
under the condition that $\E X^8$ is finite. 
The comparison study summarised below includes 
all the three choices from $(x,x^2)$ to $(x,x^2,x^3,x^4)$. 

We should also stress that other and non-polynomial choices 
are feasible and sometimes advisable. 
With $(x,\log x)$, for example, one has a specially designed gamma density. 
One also sees from (1.4)--(1.5) that the bias involved can be 
large at the tails where terms like $x^4f(x)$ can dominate over
the pilot estimate's $f''(x)$. Also, higher order moments 
need to exist to avoid large biases. 
More robust choices can easily be constructed. 

\section
{\csc 3.2. A comparison study.} 
Here we consider eight different estimators:
the kernel method $\hatt f_0$,
the second, third and fourth order \ET\ method,
the variance correction method of Jones (1991),
the multiplicative estimator of Hjort and Glad (1995), 
and finally the two local likelihood based local correctors 
on a normal start estimate given as 
$\hatt f_4$ and $\hatt f_5$ in Section 2. 
All of these have the same approximate variance, 
and bias functions of the form $\half k_2 h^2b(x)$ plus smaller order terms.
For large-sample comparison purposes it therefore 
suffices to compare the $b(x)$ functions. 

For this comparison we choose to include only the versions of 
$\hatt f_3$, $\hatt f_4$, $\hatt f_5$ that utilise 
the globally estimated normal density as start estimate, 
say $\hatt\sigma^{-1}\phi(\hatt\sigma^{-1}(x-\hatt\mu))$, 
and with the normal kernel $K=\phi$. 
Various calculations then give 
$$\eqalign{
\hatt f_3(x)&=n^{-1}\sumin h^{-1}\phi(h^{-1}(x_i-x))
        \exp[\half\{(x_i-\hatt\mu)^2-(x-\hatt\mu)^2\}/\hatt\sigma^2], \cr
\hatt f_4(x)&=\hatt f_0(x)(1+h^2/\hatt\sigma^2)^{1/2}
        \exp[-\half h^2(x-\mu)^2/\{\hatt\sigma^2(\hatt\sigma^2+h^2)\}], \cr 
\hatt f_5(x)&=\hatt f_0(x)(1+h^2/\hatt\sigma^2)^{1/2}
        \exp[-\half h^2(1+h^2/\hatt\sigma^2)
                \{\hatt f_0'(x)/\hatt f_0(x)\}^2]. \cr}$$
And by the results listed as (2.1)--(2.3) in Section 2 
the bias factor functions become 
$$\eqalign{
b_1&=f'', \cr
b_{\rm ET,2}&=f''-fg_2, \cr
b_{\rm ET,3}&=f''-fg_3, \cr
b_{\rm ET,4}&=f''-fg_4, \cr
b_2&=f''-(f')^2/f+f[\{f'/f+(x-\mu)/\sigma^2\}^2+1/\sigma^2], \cr
b_3&=f''(x)+\sigma^{-2}\{f(x)+(x-\mu)f'(x)\}, \cr
b_4&=f''(x)-f(x)\{(x-\mu)^2/\sigma^2-1\}/\sigma^2, \cr
b_5&=f''(x)-f'(x)^2/f(x)+f(x)/\sigma^2. \cr} $$
One may now display these functions for visual inspection,  
and compute useful summary numbers, for any candidate density $f$. 

For this occasion we shall be content to study the 
bias behaviour for the selection of finite normal mixtures 
put up by Marron and Wand (1992). They constructed 
15 test cases to exhibit a broad range of different density forms.
We shall limit attention to the first ten of these,
in that numbers 11--15 are too peculiar and otherworldly
for practical interest. 
The first is the standard normal, and the next nine are called  
the skewed unimodal,
the strongly skewed,
the kurtotic unimodal,
the outlier,
the bimodal,
the separated bimodal,
the skewed bimodal,
the trimodal,
and the claw. 
These mixtures are defined in terms of up to eight normal components. 
Summary quantities like $(\int b^2\,\d x)^{1/2}$, 
$\int|b|\,\d x$, $\max(|b|)$ and so on were computed. 
The study of Hjort and Jones (1996b) is broader, 
in the sense that more estimators and more summary quantities are included, 
and deeper, in the sense that also real finite-sample performance
with $h$-selectors and so on is considered. 

The table here presents just one such performance summary;
other tables were computed and gave more or less the same
main impressions. Here we give values of 
$(\int b^2)^{1/2}/\allowbreak\{\int(f'')^2\}^{1/2}$,
that is, each method's root integrated squared bias value 
normalised with that of the standard kernel method. 
To compute these for the \ET\ method one needs to 
derive formulae for the variance matrix for $(X,X^2,X^3,X^4)^\tr$ 
under a normal mixture distribution, among other quantities. 
Numerical integration was used, 
partly by necessity and partly for simplicity 
(explicit formulae for the $\int b^2$ integral may be derived 
after rather arduous calculations,
for all of the estimators except $\hatt f_5$, 
see the methods of Hjort and Glad (1996)). 
A value less than 1 is synonymous with better performance 
than the kernel estimator for (at least) all large $n$, 
as measured by the traditional mean integrated squared error criterion. 

{{\smallskip\narrower\baselineskip11pt
\parindent0pt 
{\csc Table.} \sl 
Ratios of root integrated squared bias to that of the kernel method,
for the first ten test cases of Marron and Wand (1992), 
for seven competitors; 
the second, third, fourth order Efron and Tibshirani method; 
the variance correction method of Jones;
the multiplicative kernel method of Hjort and Glad; 
and the one- and two-parametric local likelihood corrections 
on the kernel estimator. 

\smallskip
~~~~~~~~~~~~~ET-2~~~~~~ET-3~~~~~~ET-4~~~~~~~$\hatt f_2$~~~~~~~~~~$\hatt f_3$~~~~~~~~~~$\hatt f_4$~~~~~~~~~~$\hatt f_5$~~
\obeylines
\smallskip\tt
~~~1~~~~0.0000~~0.0000~~0.0000~~0.0000~~0.0000~~0.0000~~0.0000
~~~2~~~~0.8064~~0.5797~~0.5806~~0.5076~~0.4959~~0.6283~~0.5139
~~~3~~~~1.0042~~1.0117~~1.0119~~0.9966~~0.9962~~0.9993~~1.1368
~~~4~~~~0.9941~~0.9941~~0.9920~~0.9849~~0.9801~~0.9899~~1.0431
~~~5~~~~0.9947~~0.9947~~0.9891~~0.9082~~0.8844~~0.9395~~1.0340
~~~6~~~~1.0699~~1.0699~~0.9447~~0.9053~~0.8778~~0.9990~~0.8732
~~~7~~~~0.9827~~0.9827~~0.9530~~0.9336~~0.9232~~0.9914~~1.0365
~~~8~~~~1.0530~~1.0415~~1.0242~~0.9507~~0.9490~~1.0057~~0.9778
~~~9~~~~1.0219~~1.0219~~0.9941~~0.9651~~0.9580~~1.0005~~0.9599
~~10~~~~0.9974~~0.9974~~0.9974~~0.9908~~0.9895~~0.9973~~1.0188
\smallskip}}

We make the following points.

(i) The \ET\ method often wins over the kernel method, 
but it is occasionally, if the true density is too far from
normality, outperformed by its simple pilot estimator $\hatt f_0$.

(ii) It is not always advantageous to go up in polynomial order, 
not even in the asymptotic case where inclusion of such extra 
$\beta_j$s is `free' (the sampling variability of $\hatt\beta_j$s 
becomes negligible in comparison with that of $\hatt f_0$). 

(iii) Importantly, the \ET\ method is outperformed 
by some of its competitors, primarily by the $\hatt f_3$ estimator
of Hjort and Glad (1995), which wins in eight out of nine non-normal cases 
(while the normal case is a perfect draw). 
Also the quite simple variance correction
method $\hatt f_2$ of Jones (1991) fares surprisingly well 
in comparison with the more highpowered (1.3) method. 

(iv) The gains are modest once one is far from normality,
and indeed in such cases all methods considered behave quite similarly. 

(v) That such asymptotic comparisons can be a good guide for what
happens for finite and and even quite moderate $n$ is 
at least substantiated in Hjort and Glad (1996) 
for the $\hatt f_3$ estimator. 

(vi) The \ET\ method is quite flexible in that 
the user may choose $t(x)$. It remains however a 
`moment fixing procedure', and estimation methods  
$\hatt f_3$, $\hatt f_4$, $\hatt f_5$ can with greater ease
and presumably greater efficiency accommodate 
prior notions of shape for the underlying density, 
through the parametric vehicle family $f(x,\theta)$. 
In other words, the table above does not quite show 
the very best of the potential of the three last methods;
see Hjort and Jones (1996b). 

\bigskip
{\bf 4. Proof of the main result.} 
Here we prove the Proposition of Section 1. 
We begin by noting that log-likelihood function for the (1.2) model is 
$-n\log\hatt c(\beta)+\beta^\tr\sumin t(x_i)$, which means that 
$\hatt\beta$ is the minimiser of the function 
$$B_n(\beta)=\log\hatt c(\beta)-\beta^\tr\bar t,
        \quad {\rm where\ }\bar t=n^{-1}\sumin t(x_i). \eqno(4.1)$$
This is a strictly convex function, 
as long as at least $p$ distinct data points have been gathered, 
and $\hatt\beta$ is also the solution to 
$\hatt\mu_j(\beta)=\bar t_j$ for $j=1,\ldots,p$, 
where $\hatt\mu(\beta)=\E_\beta t(X)$ 
is the mean of $t(X)$ under the (1.2) model; 
$$\hatt\mu(\beta)
        =\int t(x)\hatt f_0(x)\hatt c(\beta)^{-1}\exp\{\beta^\tr t(x)\}\,\d x
        =\dell\log\hatt c(\beta)/\dell\beta. $$ 

Observe that when $n$ is moderate to large,
the start estimator $\hatt f_0$ itself will be close to the true $f$,
and the $\hatt\beta$ should with high probability not be far from zero. 
This already hints that Taylor expansions around $\beta=0$ 
will be useful for the analysis. 
One way of appreciating this is to consider the question, 
what is $\hatt\beta$ really aiming at, in this context? 
The ambition of the maximum likelihood procedure 
is to find the best parametric approximant 
to the true density in the sense of minimising 
the Kullback--Leibler distance 
$$\int f(x)\log{f(x)\over 
        \hatt f_0(x)\hatt c(\beta)^{-1}\exp\{\beta^\tr t(x)\}}\,\d x
        ={\rm const.}+\log\hatt c(\beta)-\beta^\tr \xi. $$
Let 
$$\hatt\beta_0={\rm argmin}\,A_n,
\quad {\rm where\ }A_n(\beta)=\log\hatt c(\beta)-\beta^\tr\xi. \eqno(4.2)$$
It is also defined as the solution to $\hatt\mu(\beta)=\xi$. 
Intuitively, since $\hatt f_0$ is about $O(h^2)+O_p((nh)^{-1/2})$ 
away from $f$, $\hatt\beta_0$ will be this much away from zero,
and $\hatt\beta-\hatt\beta_0$ should by familiar maximum likelihood 
theory arguments be $O_p(n^{-1/2})$. 

Both of these claims can be substantiated via the 
convexity methods of Hjort and Pollard (1996). 
The $A_n$ function converges pointwise in probability to 
$A(\beta)=\log\int f(x)\allowbreak \exp\{\beta^\tr t(x)\}\,\d x-\beta^\tr\xi$,
if this limit exists. By convexity the minimiser $\hatt\beta_0$
converges to the minimiser of the limit function, which is 0. 
This also holds in cases where $A(\beta)$ could be infinite 
for some non-null values. Under regularity conditions 
there is also convergence towards a normal limit 
for $(nh)^{1/2}(\hatt\beta_0+\half k_2h^2\gamma)$,
where $\gamma=\E t''(x)$. This can be shown by a variation
of the arguments that now follow to prove the claim about 
$\hatt\beta-\hatt\beta_0$. Study the convex function 
$$\eqalign{
n\{B_n(\hatt\beta_0+n^{-1/2}s)-B_n(\hatt\beta_0)\}
&=n\{\log\hatt c(\hatt\beta_0+n^{-1/2}s)-\log\hatt c(\hatt\beta_0)
        -n^{-1/2}s^\tr\bar t\} \cr
&=\half s^\tr\hatt\Sigma_0s-n^{1/2}s^\tr(\bar t-\xi)+r_n(s), \cr}$$
where $\hatt\Sigma_0=\dell^2\log\hatt c(0)/\dell\beta\dell\beta^\tr$
is the variance matrix of $t(X_i)$ under $\beta=0$,
and where $r_n(s)$ is seen to go to zero in probability for each $s$. 
By the results of Hjort and Pollard (1996, Sections 1 and 2),  
this is sufficient to infer that the minimiser of the function studied,
which is $n^{1/2}(\hatt\beta-\hatt\beta_0)$,
is at most $O_p(1)$ away from the minimiser of 
$\half s^\tr\hatt\Sigma_0s-n^{1/2}s^\tr(\bar t-\xi)$,
which is $\hatt\Sigma_0^{-1}n^{1/2}(\bar t-\xi)$. 
In particular $n^{1/2}(\hatt\beta-\hatt\beta_0)$ 
has a normal limit, cf.~remarks made in Section 3 in \ET\ (1996),
but in the Poisson discretisation framework. 

While useful for some purposes, exact limit distributions for 
$\hatt\beta$ and $\hatt\beta_0$ do not directly concern us here.
But the above gives us Taylor approximations for necessary 
quantities around $\beta=0$, which turn out to be good enough,
in that $\hatt\beta$ and $\hatt\beta_0$ both go to zero when
$h\arr0$ and $nh\arr\infty$. 
In addition to $\hatt\Sigma_0$ we shall need to work 
with the mean of $t(X)$ under $\beta=0$, that is, 
under the $\hatt f_0(x)$ density. These are easy to approximate,
$$\eqalign{
\hatt\xi_0=\E_0t(X)&=n^{-1}\sumin \int t(x)K_h(x_i-x)\,\d x \cr
&=n^{-1}\sumin \int K(u)t(x_i+hu)\,\d u 
        =\bar t+\half k_2h^2\bar{t''}+O_p(h^4), \cr} \eqno(4.3)$$
where $\bar{t''}$ is the $p$-vector with 
components $n^{-1}\sumin t_j''(x_i)$. Similarly, 
$$\hatt\Sigma_0=\Var_0t(X)
 =n^{-1}\sumin\int K(u)(tt^\tr)(x_i+hu)\,\d u-\hatt\xi_0\hatt\xi_0^\tr 
 =\hatt\Sigma+\half k_2h^2\hatt C+O_p(h^4), \eqno(4.4)$$
say, where $\hatt\Sigma$ is the empirical variance matrix for 
$t(x_1),\ldots,t(x_n)$, and with a certain expression available 
for the $\hatt C$ matrix from Taylor expansion. 

We are now in a position to approximate $\hatt c(\beta)$
and hence $A_n(\beta)$ and $B_n(\beta)$ and their minima. One finds 
$$\hatt c(\beta)\doteq 1+\beta^\tr\hatt\xi_0
+\half\beta^\tr(\hatt\Sigma_0+\hatt\xi_0\hatt\xi_0^\tr)\beta
        +\sixth\int \{\beta^\tr t(x)\}^3\hatt f_0(x)\,\d x, $$
which implies 
$B_n(\beta)=\beta^\tr(\hatt\xi_0-\bar t)
        +\half\beta^\tr\hatt\Sigma_0\beta+O_p(\|\beta\|^3)$. 
Hence, to a first order of approximation, 
$$\hatt\beta\doteq \hatt\Sigma_0^{-1}(\bar t-\hatt\xi_0)
        \doteq -\half k_2h^2\hatt\Sigma^{-1}\bar{t''}. $$
This yields a little string of approximations for the
\ET\ estimator, useful for different purposes: 
$$\eqalign{\hatt f(x)
&\doteq \hatt f_0(x)
{1+\hatt\beta^\tr t(x)+\half\hatt\beta^\tr t(x)t(x)^\tr\hatt\beta
\over 1+\hatt\beta^\tr\hatt\xi_0
        +\half\hatt\beta^\tr(\hatt\Sigma_0
        +\hatt\xi_0\hatt\xi_0^\tr)\hatt\beta} \cr 
&\doteq \hatt f_0(x)\{1+\hatt\beta^\tr(t(x)-\hatt\xi_0)
        +\half\hatt\beta^\tr(t(x)-\hatt\xi_0)(t(x)-\hatt\xi_0)'
                -\hatt\Sigma_0\}\hatt\beta\} \cr 
&\doteq \hatt f_0(x)\{1-\half k_2h^2(\bar{t''})'\hatt\Sigma^{-1}
        (t(x)-\bar t)\}. \cr} \eqno(4.5)$$ 
Essentially this shows that the correction term 
$\hatt c(\hatt\beta)^{-1}\exp\{\hatt\beta^\tr t(x)\}$ 
in the (1.3) estimator can be represented as 
$1-\half k_2h^2\{O_p(h^2)+\hatt g(x)\}$, where $\hatt g(x)$ is 
an ordinary $n^{1/2}$-consistent estimator of the $g(x)$ function (1.4),
when $h\arr0$, and the $O_p(h^2)$ term stemming from (4.3) and (4.4). 
Inspecting matters involved in (4.3)--(4.5) 
shows that $\hatt f$ admits a representation of the form 
$$\hatt f(x)=\hatt f_0(x)
        [1-\half k_2h^2\{g(x)+n^{-1}d(x)+O_p(h^2+n^{-2})
                +n^{-1/2}\tau(x)V_n(x)\}], \eqno(4.6)$$
say, where $n^{-1}d(x)$ and $n^{-1}\tau(x)^2$ are 
the dominant terms in the bias and variance of $\hatt g(x)$
as an estimator of $g(x)$, 
and where $V_n(x)$ has mean zero and variance going to 1. 

This is quickly seen to imply the main result (1.5), in that 
$\hatt f_0=f+\half k_2h^2f''+o(h^2)+(nh)^{-1/2}R(K)^{1/2}f(x)^{1/2}U_n$,
in which $U_n(x)$ has mean zero and variance going towards 1.  

\bigskip
{\bf 5. Supplementing remarks.} 
This final section offers some extraneous results,
generalisations and remarks,
some of which point to possible further research. 

\section
{\csc 5.1. Influence of $\beta$ estimates.}
The necessary requirements for consistency of 
the (1.3) estimator are $h\arr0$ and $nh\arr\infty$,
and under these conditions result (1.5) 
satisfactorily explains its large-sample behaviour. 
One feature of this result is that the sampling variability
of $(\hatt\beta_1,\ldots,\hatt\beta_p)^\tr$ 
only has a minor, secondary effect on the density estimator. 
This is not the complete story, however, 
and Efron and Tibshirani's treatment sheds more light on
finite-sample behaviour, and particular provides diagnostic tools
to monitor effects of different $t(x)$ choices. 
One should note that if too many $\beta_j$s 
are allowed in the model, then the variance 
of $\hatt c(\hatt\beta)^{-1}\exp\{\hatt\beta^\tr t(x)\}$ 
will be more influential than what the large-sample result suggests. 
The point is partly that the $\tau(x)$ term in (4.6) becomes large. 

\section
{\csc 5.2. Choosing the bandwidth.} 
There is some ongoing cross-fertilisation regarding the 
issue of bandwidth selection in density estimation and regression areas,
traces of which can be seen in Efron and Tibshirani (1996), 
Fan and Gijbels (1996), Simonoff (1996), Loader (1996b) 
and Hjort and Jones (1996b). 
This is also important in view of the many new proposals 
for density estimators, where the literature is bewilderingly overcrowded 
with approaches in the very simplest case of the kernel method,
but rather silent regarding other methods. 
Recent and partly conflicting advice can be found in 
for example Jones, Marron and Sheather (1995) and Loader (1996b). 

The present point is to indicate how the density estimation viewpoint,
and result (1.5) in particular, can be utilised for the \ET\ method,
and indeed also for the other estimators presented in Section 2. 
Least squares cross validation remains a generally valid option, 
even though it loses the `unbiasedness' property it enjoys 
(see for example Scott, 1992, or Wand and Jones, 1995)
in the case of the kernel method. 
Another method, which typically would display less sampling variability, 
is based on minimising the approximate 
mean integrated (weighted) squared error,
which in terms of a method with $\half k_2h^2b(x)$ bias is  
$\quart k_2^2h^4\int wb^2\,\d x+(nh)^{-1}R(K)\int wf\,\d x$, 
with a suitable weight function $w$. 
This could the constant 1, giving the familiar `mise' 
and then a global bandwidth, 
or it could be a window function around a given $x$, 
leading to a local bandwidth. The minimiser is 
$$h_0=\Bigl\{{R(K)\over k_2^2}
        {\int wf\,\d x\over \int wb^2\,\d x}\Bigr\}^{1/5}\,n^{-1/5}
     =\Bigl({\int wb_0^2\,\d x\over \int wb^2\,\d x}\Bigr)^{1/5}
        \,h_{0,{\rm trad}}, $$
say, where $b_0=f''$ for the traditional kernel method. 
Thus one could think of modifying the best bandwidth for the kernel method 
with the factor $(\int wb_0^2/\int wb^2)^{1/5}$. 
In view of the Table of Section 3.2 this factor will not be very far from 1,
so one cannot go very wrong by using a global or local bandwidth 
selected for the kernel method. More sophisticatedly one might 
estimate $\int wb^2=\int w(f''-fg)^2$ using a pilot estimate and 
deduction of the positive bias involved. 

We also point out that the minimum possible value of the 
criterion function here can be expressed as 
$k_2^{2/5}R(K)^{4/5}/n^{4/5}$ 
times an expression in $f$ alone. 
The best kernel function to use, in this sense, 
is the $K$ which minimises $k_2^{2/5}R(K)^{4/5}$.
The solution is the Yepanechnikov kernel 
${3\over 2}\{1-(2x)^2\}\,I\{|x|\le\half\}$, 
or any scaled version thereof. This result is valid 
for all of the estimators considered in Section 2 and 3. 

\section
{\csc 5.3. General start estimator.} 
The (1.2)--(1.3) recipe will of course work with other start estimators 
than the kernel method. Suppose that a $\hatt f_0$ 
different from the kernel estimator is used, 
say one of the competitors listed in Section 2,
with approximate bias $\half k_2h^2b(x)$ 
and approximate integral $1+dh^2$, where $d=\int b(x)\,\d x$. 
This leads to new versions of 
$\hatt c(\beta)$, $\hatt\xi_0$, $\hatt\Sigma_0$ and so on,
defined analogously to those in Section 4 for the kernel method. 
For the analogue of the (4.1) function the arguments of Section 4 give 
$$B_n(\beta)=\log\hatt c(\beta)-\beta^\tr\bar t
        =\beta^\tr\{(1-dh^2)\hatt\xi_0-\bar t\}
         +\half\beta^\tr\hatt\Sigma_0\beta+d h^2+O_p(h^4+\|\beta\|^3), $$
with minimiser 
$\hatt\beta\doteq\hatt\Sigma_0^{-1}\{\bar t-(1-dh^2)\hatt\xi_0\}$. 
Under appropriate regularity conditions, which should then be checked
for the candidate $\hatt f_0$ in question, one would have 
$$\hatt c(\hatt\beta)^{-1}\exp\{\hatt\beta^\tr t(x)\}
        \doteq 1-\half k_2h^2\{\gamma(b)^\tr\hatt\Sigma^{-1}(t(x)-\bar t)
                +O_p(h^2+(nh)^{-1})\}, $$
where $\gamma(b)=\int b(x)t(x)\,\d x-(2d/k_2)\xi$,
by the reasoning of Section 4 again. This leads to 
an approximate bias of the form $\half k_2h^2b_{\rm new}(x)$
for this generalised \ET\ estimator, where 
$$b_{\rm new}=b-f\gamma(b)^\tr\Sigma^{-1}(t(x)-\xi), $$
and the approximate variance is again 
approximately unchanged, as with (1.5). 

The previous case corresponds to $b=f''$, for which $d=0$
and $\gamma(b)=\int f''t\,\d x=\int ft''\,\d x$ by partial integration. 
We note in particular that if the \ET\ method is iterated once,
then $b=f''-fg$ gives, after some analysis, $b_{\rm new}=b$ unchanged.
Inspection reveals that such an iteration simply gives 
$\hatt\beta=0$ at the second stage, so the formal iteration does not
change the estimator in this case. For other choices of $\hatt f_0$,
like the multiplicative $\hatt f_3$ which appeared to win the
competition in Section 3, this scheme gives new estimators 
with new bias functions. 

\section
{\csc 5.4. Similar ideas for regression.}
Ideas and methodology influence each other in the related worlds
of density and regression estimation, as exemplified by
the \ET\ article. So there ought to be a parallel 
to the (1.3) method for nonparametric regression.
One version is as follows. 
Take any reasonable $\hatt m_0$ smooth, 
say a local linear or local quadratic kernel estimator, 
see for example Wand and Jones (1995) or Fan and Gijbels (1996).
Then multiply with $\hatt\beta_0+\hatt\beta_1x+\hatt\beta_2x^2$, 
say, with these coefficients arrived at via direct least squares for 
$y_i-\hatt m_0(x_i)(\beta_0+\beta_1x_i+\beta_2x_i^2)$,
which amounts to weighted least squares for 
$y_i/\hatt m_0(x_i)$ on $\beta_0+\beta_1x_i+\beta_2x_i^2$. 
Some caution would be called for in cases where 
$\hatt m_0$ crosses zero. This method would complement Glad's (1995) method;
she starts with a parametric estimator 
and multiplies with a nonparametric kernel type correction. 
The aim is to reduce bias while keeping the variance level. 

\section
{\csc 5.5. Higher dimension.}
Estimation methods studied in this paper 
and results on their behaviour can be generalised 
to vector data without too much difficulties.
In particular result (1.5) can be generalised to 
the situation discussed in Efron and Tibshirani's Section 4. 
Presumably the potential of these semiparametric methods 
as significant improvers on standard methods is larger 
in higher dimensions, where purely nonparametric estimators
are in serious trouble. 

\section
{\csc 5.6. Generalisations.} 
There are still further variations on the theme. 
One possibility worth mentioning is to use the (1.2)--(1.3) scheme
with a parametric start estimator, again to `fix moments'.
For such use it would be best to truncate the range of the 
start estimator, to avoid trouble with infinite integrals. 

Finally we outline an approach that may allow many terms in 
$(t_1(x),\ldots,t_p(x))^\tr$,
to allow for greater flexibility in the functional form,
but restraining the coefficients by suitable penalisation,
to avoid too much noise. This can be afforded 
by maximising the penalised log-likelihood function 
$$-n\log\Bigl[\int\hatt f_0(x)
        \exp\Bigl\{\sum_{j=1}^p\beta_jt_j(x)\Bigr\}\,\d x\Bigr]
        +n\sum_{j=1}^p\beta_j\bar t_j
        -\half\sum_{j=1}^p\beta_j^2/\tau_j^2, $$
where the $\tau_j$s should go to zero with increasing $j$,
say as $\tau_j=\tau_0/j$. This would correspond to 
a Bayesian prior which places independent normal $(0,\tau_j^2)$ distributions
on the $\beta_j$s. This is in fact one of several approaches 
towards nonparametric Bayesian density estimation 
laid out in Hjort (1995).
Section 5 of that paper gives further details and also a 
Schwarz type information criterion to help decide 
on the number of terms used. 



\bigskip
{\bf Acknowledgements.}
In connection with this work I have benefited 
from comments from Ingrid Glad and Jeff Simonoff, 
and in particular from Chris Jones. 


\bigskip
\parindent0pt
\parskip3pt
\baselineskip12pt

\centerline{\bf References}

\medskip
\ref{%
Efron, B.~and Tibshirani, R. (1996).
Using specially designed exponential families for density estimation.
{\sl Annals of Statistics} {\bf 24}, to appear.} 

\ref{%
Fan, J.~and Gijbels, I. (1996).
{\sl Local Polynomial Modelling and its Applications.}
Chapman and Hall, London. To exist.} 


\ref{%
Glad, I.K. (1995). 
Parametrically guided nonparametric regression.
Department of Mathematical Sciences, 
Norwegian Institute of Technology.} 

\ref{%
Hjort, N.L. (1995). 
Bayesian approaches to semiparametric density estimation.
In {\sl Bayes\-ian Statistics 5}, 
proceedings of the Fifth International Valencia Meeting on Bayesian Statistics 
(eds.~J.~Berger, J.Bernardo, A.P.~Dawid, A.F.M.~Smith),
xxx--xxx. Oxford University Press.}

\ref{%
Hjort, N.L. (1996).
Multiplicative higher order bias kernel density estimators.
To exist.} 

\ref{%
Hjort, N.L.~and Glad, I.K. (1995).
Nonparametric density estimation with a parametric start.
{\sl Annals of Statistics} {\bf 23}, 882--904.}

\ref{%
Hjort, N.L.~and Glad, I.K. (1996).
Exact performance for a multiplicative semiparametric density estimator, 
for normal mixture truths. To exist.}  

\ref{%
Hjort, N.L.~and Jones, M.C. (1996a).
Locally parametric nonparametric density estimation.
{\sl Annals of Statistics} {\bf 24}, to appear.} 

\ref{%
Hjort, N.L.~and Jones, M.C. (1996b). 
Recent advances in semiparametric kernel density estimation: 
A review and comparison. 
To exist.} 
 
\ref{%
Hjort, N.L.~and Pollard, D. (1996).
Asymptotics for minimisers of convex processes. 
{\sl Annals of Statistics}, to appear.}

\ref{%
Jones, M.C. (1991).
On correcting for variance inflation in kernel density estimation.
{\sl  Computational Statistics \& Data Analysis} {\bf 11}, 3--15.}

\ref{%
Jones, M.C., Linton, O.~and Nielsen, J.P. (1995).
A simple bias reduction method 
for density and regression estimation.
{\sl Biometrika} {\bf 82}, 327--338.}  

\ref{%
Jones, M.C., Marron, J.S.~and Sheather, S.J. (1995).
A brief survey of bandwidth selection.
{\sl Journal of the American Statistical Association} {\bf 90}, 
to appear.} 

\ref{%
Loader, C.R. (1996a). 
Local likelihood density estimation. 
{\sl Annals of Statistics}, to appear.} 

\ref{%
Loader, C.R. (1996b). 
Old Faithful erupts: Bandwidth selection reviewed.
Manuscript, AT\&T Laboratories, Murray Hill.}

\ref{%
Marron, J.S.~and Wand, M.P. (1992). 
Exact mean integrated squared error. 
{\sl Annals of Statistics} {\bf 20}, 712--736.}

\ref{%
Peirce, B. (1873). 
On the theory of errors of observations.
Appendix No.~21 (pp.~200--224 and plate 27) of 
{\sl Report of the Superintendent of the U.S.~Coast Survey 
for the year ending June 1870}. G.P.O., Washington.} 


\ref{%
Scott, D.W. (1992).
{\sl Multivariate Density Estimation:
Theory, Practice, and Visualization.}
Wiley, New York.}

\ref{%
Simonoff, J. (1996).
{\sl Smoothing in Statistics.}
Springer, to exist.}

\ref{%
Wand, M.P.~and Jones, M.C. (1995).
{\sl Kernel Smoothing.}
Chapman \& Hall, London.}

\bigskip
\bigskip
\parindent20pt
\baselineskip13pt

\bye